\documentclass[reqno]{amsart}
\usepackage{setspace}
\usepackage{graphicx}
\usepackage{amsmath}
\usepackage{amsfonts}
\usepackage{amssymb}

\begin{document}

\centerline{\huge \bf Multiple Lattice Tilings in Euclidean Spaces}

\bigskip\medskip
\centerline{Dedicated to Professor Dr. Christian Buchta on the occasion of his 60th birthday}

\bigskip\medskip
\centerline{\large\bf Qi Yang and Chuanming Zong}

\vspace{1cm}
\centerline{\begin{minipage}{12.8cm}
{\bf Abstract.} In 1885, Fedorov discovered that a convex domain can form a lattice tiling of the Euclidean plane if and only if it is a parallelogram or a centrally symmetric hexagon. This paper proves the following results: {\it Besides parallelograms and centrally symmetric hexagons, there is no other convex domain which can form a two-, three- or four-fold lattice tiling in the Euclidean plane. However, there are both octagons and decagons which can form five-fold lattice tilings. Whenever $n\ge 3$, there are non-parallelohedral polytopes which can form five-fold lattice tilings in the $n$-dimensional Euclidean space.}
\end{minipage}}


\vspace{1cm}
\noindent
{\Large\bf 1. Introduction}

\bigskip\noindent
Planar tilings is an ancient subject in our civilization. It has been considered in the arts by craftsmen since antiquity. Up to now, it is still an active research field in mathematics and some basic problems remain unsolved. In 1885, Fedorov \cite{fedo} discovered that there are only two types of two-dimensional lattice tiles: {\it parallelograms and centrally symmetric hexagons}. In 1917, Bieberbach suggested Reinhardt (see \cite{rein}) to determine all the two-dimensional convex congruent tiles. However, to complete the list turns out to be challenging and dramatic. Over the years, the list has been successively extended by Reinhardt, Kershner, James, Rice, Stein, Mann, McLoud-Mann and Von Derau (see \cite{zong14,mann}), its completeness has been mistakenly announced several times! In 2017, M. Rao \cite{rao} announced a completeness proof based on computer checks.

The three-dimensional case was also studied in the ancient time. More than 2,300 years ago, Aristotle claimed that both identical regular tetrahedra and identical cubes can fill the whole space without gap. The cube case is obvious! However, the tetrahedron case is wrong and such a tiling is impossible (see \cite{lazo}).

Let $K$ be a convex body with (relative) interior ${\rm int}(K)$, (relative) boundary $\partial (K)$ and volume ${\rm vol}(K)$, and let $X$ be a discrete set, both in $\mathbb{E}^n$. We call $K+X$ a {\it translative tiling} of $\mathbb{E}^n$ and call $K$ a {\it translative tile} if $K+X=\mathbb{E}^n$ and the translates ${\rm int}(K)+{\bf x}_i$ are pairwise disjoint. In other words, if $K+X$ is both a packing and a covering in $\mathbb{E}^n$. In particular, we call $K+\Lambda$ a {\it lattice tiling} of $\mathbb{E}^n$ and call $K$ a {\it lattice tile} if $\Lambda $ is an $n$-dimensional lattice. Apparently, a translative tile must be a convex polytope. Usually, a lattice tile is called a {\it parallelohedron}.

In 1885, Fedorov \cite{fedo} also characterized the three-dimensional lattice tiles: {\it A three-dimensional lattice tile must be a parallelotope, an hexagonal prism, a rhombic dodecahedron, an elongated dodecahedron, or a truncated octahedron.} The situations in higher dimensions turn out to be very complicated. Through the works of Delone \cite{delo}, $\check{S}$togrin \cite{stog} and Engel \cite{enge}, we know that there are exact $52$ combinatorially different types of parallelohedra in $\mathbb{E}^4$. A computer classification for the five-dimensional parallelohedra was announced by Dutour Sikiri$\acute{\rm c}$, Garber, Sch$\ddot{\rm u}$rmann and Waldmann \cite{dgsw} only in 2015.

Let $\Lambda $ be an $n$-dimensional lattice. The {\it Dirichlet-Voronoi cell} of $\Lambda $ is defined by
$$C=\left\{ {\bf x}: {\bf x}\in \mathbb{E}^n,\ \| {\bf x}, {\bf o}\|\le \| {\bf x}, \Lambda \|\right\},$$
where $\| X, Y\|$ denotes the Euclidean distance between $X$ and $Y$. Clearly, $C+\Lambda $ is a lattice tiling and the Dirichlet-Voronoi cell $C$ is a parallelohedron. In 1908,
Voronoi \cite{voro} made a conjecture that {\it every parallelohedron is a linear transformation image of the Dirichlet-Voronoi cell of a suitable lattice.} In $\mathbb{E}^2$, $\mathbb{E}^3$ and $\mathbb{E}^4$, this conjecture was confirmed by Delone \cite{delo} in 1929. In higher dimensions, it is still open.

To characterize the translative tiles is another fascinating problem. First it was shown by Minkowski \cite{mink} in 1897 that {\it every translative tile must be centrally symmetric}. In 1954, Venkov \cite{venk} proved that {\it every translative tile must be a lattice tile $($parallelohedron$)$} (see \cite{alek} for generalizations). Later, a new proof for this beautiful result was independently discovered by McMullen \cite{mcmu}.

Let $X$ be a discrete multiset in $\mathbb{E}^n$ and let $k$ be a positive integer. We call $K+X$ a {\it $k$-fold translative tiling} of $\mathbb{E}^n$ and call $K$ a {\it $k$-fold translative tile} if every point ${\bf x}\in \mathbb{E}^n$ belongs to at least $k$ translates of $K$ in $K+X$ and every point ${\bf x}\in \mathbb{E}^n$ belongs to at most $k$ translates of ${\rm int}(K)$ in ${\rm int}(K)+X$. In other words, $K+X$ is both a $k$-fold packing and a $k$-fold covering in $\mathbb{E}^n$. In particular, we call $K+\Lambda$ a {$k$-fold lattice tiling} of $\mathbb{E}^n$ and call $K$ a {\it $k$-fold lattice tile} if $\Lambda $ is an $n$-dimensional lattice. Apparently, a $k$-fold translative tile must be a convex polytope. In fact, as it was shown by Gravin, Robins and Shiryaev \cite{grs}, {\it a $k$-fold translative tile must be a centrally symmetric polytope with centrally symmetric facets.}

Multiple tilings were first investigated by Furtw\"angler \cite{furt} in 1936 as a generalization of Minkowski's conjecture on cube tilings. Let $C$ denote the $n$-dimensional unit cube. Furtw\"angler made a conjecture that {\it every $k$-fold lattice tiling $C+\Lambda$ has twin cubes. In other words, every multiple lattice tiling $C+\Lambda$ has two cubes sharing a whole facet.} In the same paper, he proved the two- and three-dimensional cases. Unfortunately, when $n\ge 4$, this beautiful conjecture was disproved by Haj\'os \cite{hajo} in 1941. In 1979, Robinson \cite{robi} determined all the integer pairs $\{ n,k\}$ for which Furtw\"angler's conjecture is false. We refer to Zong \cite{zong05,zong06} for an introduction account and a detailed account on this fascinating problem, respectively, to pages 82-84 of Gruber and Lekkerkerker \cite{grub} for some generalizations.

\medskip
In 1994, Bolle \cite{boll} proved that {\it every centrally symmetric lattice polygon is a multiple lattice tile}. Let $\Lambda $ denote the two-dimensional integer lattice, and let $D_8$ denote the octagon with vertices $(1,0)$, $(2,0)$, $(3,1)$, $(3,2)$, $(2,3)$, $(1,3)$, $(0,2)$ and $(0,1)$. As a particular example of Bolle's theorem, it was discovered by Gravin, Robins and Shiryaev \cite{grs} that {\it $D_8+\Lambda$ is a seven-fold lattice tiling of $\mathbb{E}^2$.} Apparently, the octagon $D_8$ is not a lattice tile. Based on this example and McMullen's criterion on parallelohedra (see Lemma 3 in Section 3), one can easily deduce that, whenever $n\ge 2$, there is a non-parallelohedral polytope which can form a seven-fold lattice tiling in $\mathbb{E}^n$.

\medskip
In 2000, Kolountzakis \cite{kolo} proved that, if $D$ is a two-dimensional convex domain which is not a parallelogram and $D+X$ is a multiple tiling in $\mathbb{E}^2$, then $X$ must be a finite union of translated two-dimensional lattices. In 2013, a similar result in $\mathbb{E}^3$ was discovered by Gravin, Kolountzakis, Robins and Shiryaev \cite{gkrs}.

Let $P$ denote an $n$-dimensional centrally symmetric convex polytope, let $\tau (P)$ denote the smallest integer $k$ such that $P$ is a $k$-fold translative tile, and let $\tau^* (P)$ denote the smallest integer $k$ such that $P$ is a $k$-fold lattice tile. For convenience, we define $\tau (P)=\infty$ (or $\tau^*(P)=\infty$) if $P$ can not form translative tiling (or lattice tiling) of any multiplicity. Clearly, for every centrally symmetric convex polytope we have
$$\tau (P)\le \tau^*(P).$$

\medskip
It is a basic and natural problem to study the distribution of the integers $\tau (P)$ or $\tau^*(P)$, when $P$ runs over all $n$-dimensional polytopes. In particular, is there an $n$-dimensional polytope $P$ satisfying $\tau (P)=2$ or $3$? is there an $n$-dimensional polytope $P$ satisfying $\tau (P)\not= \tau^*(P)$? is there a convex domain $D$ satisfying $2\le \tau^*(D)\le 6$?

\medskip
In this paper, we will prove the following results.

\medskip
\noindent
{\bf Theorem 1.} {\it If $D$ is a two-dimensional centrally symmetric convex domain which is neither a parallelogram nor a centrally symmetric hexagon, then we have
$$\tau^*(D)\ge 5,$$
where the equality holds when $D$ is a suitable octagon or a suitable decagon.}

\medskip\noindent
{\bf Theorem 2.} {\it Whenever $n\ge 3$, there is a non-parallelohedral polytope $P$ such that $P+\mathbb{Z}^n$ is a five-fold lattice tiling of the $n$-dimensional Euclidean space.}

\vspace{0.6cm}
\noindent
{\Large\bf 2. A General Lower Bound and Two Particular Examples}

\bigskip\noindent
In this section, we will prove Theorem 1. First, let's recall some basic results which will be useful in this paper.

Let $D$ denote a two-dimensional centrally symmetric convex domain and let $\sigma $ be a non-singular affine linear transformation from $\mathbb{E}^2$ to $\mathbb{E}^2$. If $D+X$ is a $k$-fold translative tiling of $\mathbb{E}^2$, then $\sigma (D)+\sigma (X)$ is a $k$-fold translative tiling of $\mathbb{E}^2$ as well. Therefore, we have
$$\tau (D)=\tau (\sigma (D))$$
and
$$\tau^* (D)=\tau^* (\sigma (D)).$$
In particular, without loss of generality, we may assume $\Lambda =\mathbb{Z}^2$ when we study multiple lattice tilings $D+\Lambda $.

\medskip
Let $\delta_k(D)$ denote the density of the densest $k$-fold lattice packings of $D$ and let $\theta_k(D)$ denote the density of the thinnest $k$-fold lattice coverings of $D$. Clearly, $\delta_1(D)$ is the density $\delta (D)$ of the densest lattice packings of $D$ and $\theta_1(D)$ is the density $\theta (D)$ of the thinnest lattice coverings of $D$. Dumir and Hans-Gill \cite{dumi} and G. Fejes T\'oth \cite{feje} proved the following result.

\medskip\noindent
{\bf Lemma 1.} {\it If $k=2,$ $3$ or $4$, then
$$\delta_k(D)=k\cdot \delta (D)$$
holds for every two-dimensional centrally symmetric convex domain $D$.}

\medskip
In 1994, Bolle \cite{boll} proved the following criterion for the two-dimensional multiple lattice tilings.

\medskip\noindent
{\bf Lemma 2.} {\it A convex polygon is a $k$-fold lattice tile for a lattice $\Lambda$ and some positive integer $k$ if and only if the following conditions are satisfied:

\noindent
{\bf 1.} It is centrally symmetric.

\noindent
{\bf 2.} When it is centered at the origin, in the relative interior of each edge $G$ there is a point of ${1\over 2}\Lambda $.

\noindent
{\bf 3.} If the midpoint of $G$ is not in ${1\over 2}\Lambda $, then $G$ is a lattice vector of $\Lambda $.}

\medskip
Let $P_{2m}$ denote a centrally symmetric convex $2m$-gon centered at the origin, let $\mathcal{P}_{2m}$ denote the set of all such $2m$-gons, and let $\mathcal{D}$ denote the family of all two-dimensional centrally symmetric convex domains. It follows by Fedorov and Venkov's results that
$$\tau (D)=\tau^*(D)=1$$
if and only if $D\in \mathcal{P}_4\cup \mathcal{P}_6$. Then, it would be both important and interesting to determine the values of
$$\min_{D\in \mathcal{D}\setminus \{\mathcal{P}_4\cup \mathcal{P}_6\}}\tau (D)$$
and
$$\min_{D\in \mathcal{D}\setminus \{\mathcal{P}_4\cup \mathcal{P}_6\}}\tau^* (D).\eqno (1)$$
Theorem 1 determines the value of (1).

\bigskip\noindent
{\bf Proof of Theorem 1.} Let $k$ be a positive integer satisfying $k\le 4$. If $D$ is a two-dimensional centrally symmetric convex domain which can form a $k$-fold lattice tiling in the Euclidean plane, then we have
$$\delta_k(D)=k.$$
By lemma 1 it follows that
$$\delta (D)={\delta_k(D)\over k}=1$$
and therefore $D$ must be a parallelogram or a centrally symmetric hexagon. In other words, if $D$ is neither a parallelogram nor a centrally symmetric hexagon, then we have
$$\tau^*(D)\ge 5.\eqno(2)$$

We take $\Lambda =\mathbb{Z}^2$. As shown by Figure 1, let $D_8$ denote the octagon with vertices
$$\begin{array}{ll}
{\bf v}_1=\left(-\mbox{${3\over {10}}$}, -2\right),\hspace{1cm} & {\bf v}_2=\left(\mbox{${3\over {10}}$}, -1\right),\\
\vspace{-0.3cm}
& \\
{\bf v}_3=\left(\mbox{${7\over {10}}$}, 0\right), & {\bf v}_4=\left(\mbox{${{13}\over {10}}$}, 2\right),\\
\vspace{-0.3cm}
&\\
{\bf v}_5=\left(\mbox{${3\over {10}}$}, 2\right), & {\bf v}_6=\left(-\mbox{${3\over {10}}$}, 1\right),\\
\vspace{-0.3cm}
&\\
{\bf v}_7=\left(-\mbox{${7\over {10}}$}, 0\right), & {\bf v}_8=\left(-\mbox{${{13}\over {10}}$}, -2\right).
\end{array}$$

\begin{figure}[!ht]
\centering
\includegraphics[scale=0.55]{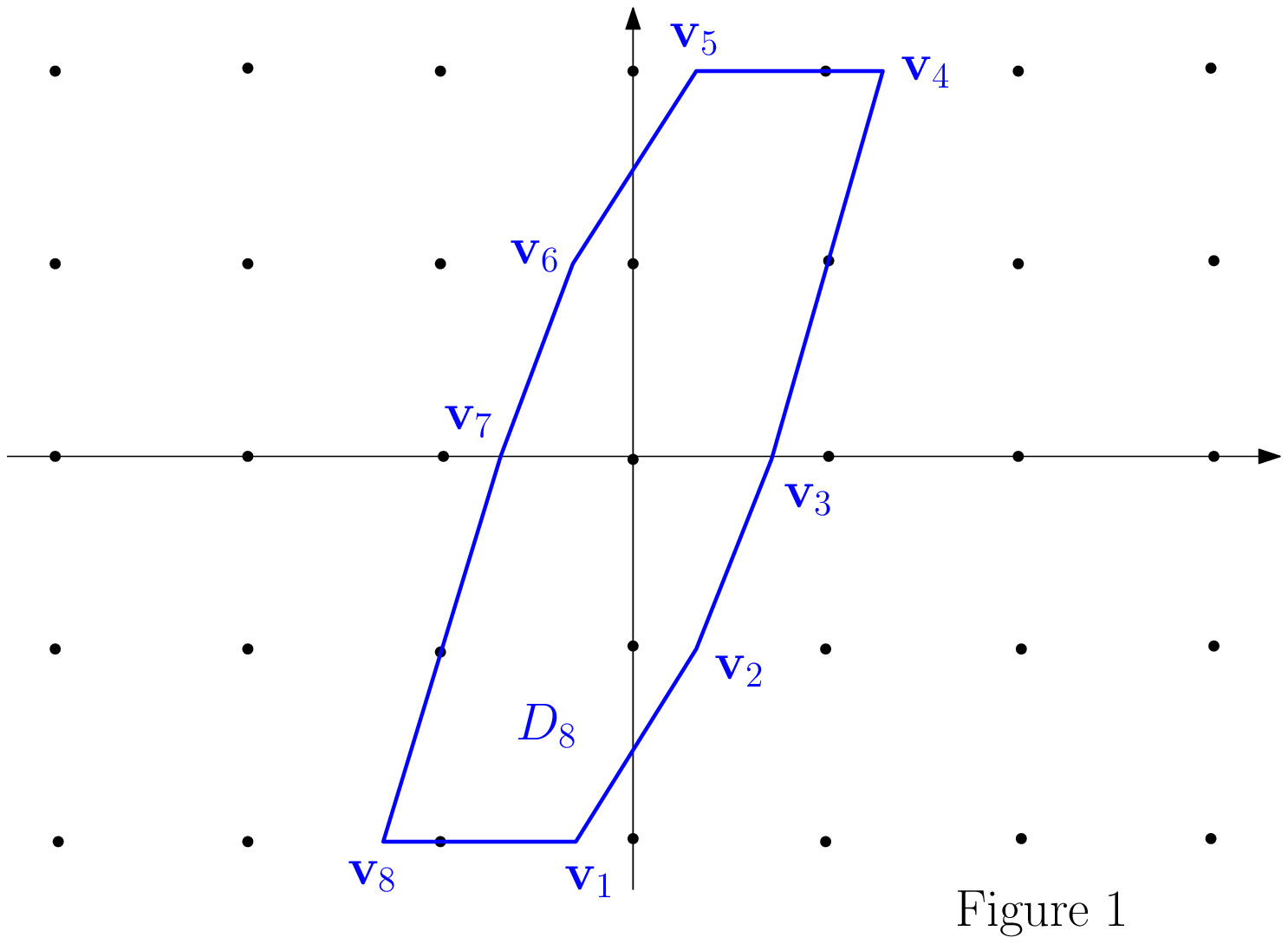}
\end{figure}
It can be easily verified that
$${\bf u}_i=\mbox{${1\over 2}$}({\bf v}_i+{\bf v}_{i+1})\in \mbox{${1\over 2}$}\Lambda, \quad i=1, 2, \ldots, 8, $$
where ${\bf v}_9={\bf v}_1$, and
$${\rm vol}(D_8)=5.$$
It follows from Lemma 2 that $D_8+\Lambda$ is a five-fold lattice tiling. Combined with (2), it can be deduced that
$$\tau^*(D_8)=5.$$

Similarly, as shown by Figure 2, let $D_{10}$ denote the decagon with vertices
$$\begin{array}{ll}
{\bf v}_1=\left(-\mbox{${3\over 5}$}, -\mbox{${5\over 4}$}\right),\hspace{1cm} & {\bf v}_2=\left(\mbox{${3\over 5}$}, -\mbox{${3\over 4}$}\right),\\
\vspace{-0.3cm}
& \\
{\bf v}_3=\left(\mbox{${7\over 5}$}, -\mbox{${1\over 4}$}\right), & {\bf v}_4=\left(\mbox{${8\over 5}$}, \mbox{${1\over 4}$}\right),\\
\vspace{-0.3cm}
&\\
{\bf v}_5=\left(\mbox{${7\over 5}$}, \mbox{${3\over 4}$}\right), & {\bf v}_6=\left(\mbox{${3\over 5}$}, \mbox{${5\over 4}$}\right),\\
\vspace{-0.3cm}
&\\
{\bf v}_7=\left(-\mbox{${3\over 5}$}, \mbox{${3\over 4}$}\right), & {\bf v}_8=\left(-\mbox{${7\over 5}$}, \mbox{${1\over 4}$}\right),\\
\vspace{-0.3cm}
&\\
{\bf v}_9=\left(-\mbox{${8\over 5}$}, -\mbox{${1\over 4}$}\right), & {\bf v}_{10}=\left(-\mbox{${7\over 5}$}, -\mbox{${3\over 4}$}\right).
\end{array}$$

\begin{figure}[!ht]
\centering
\includegraphics[scale=0.5]{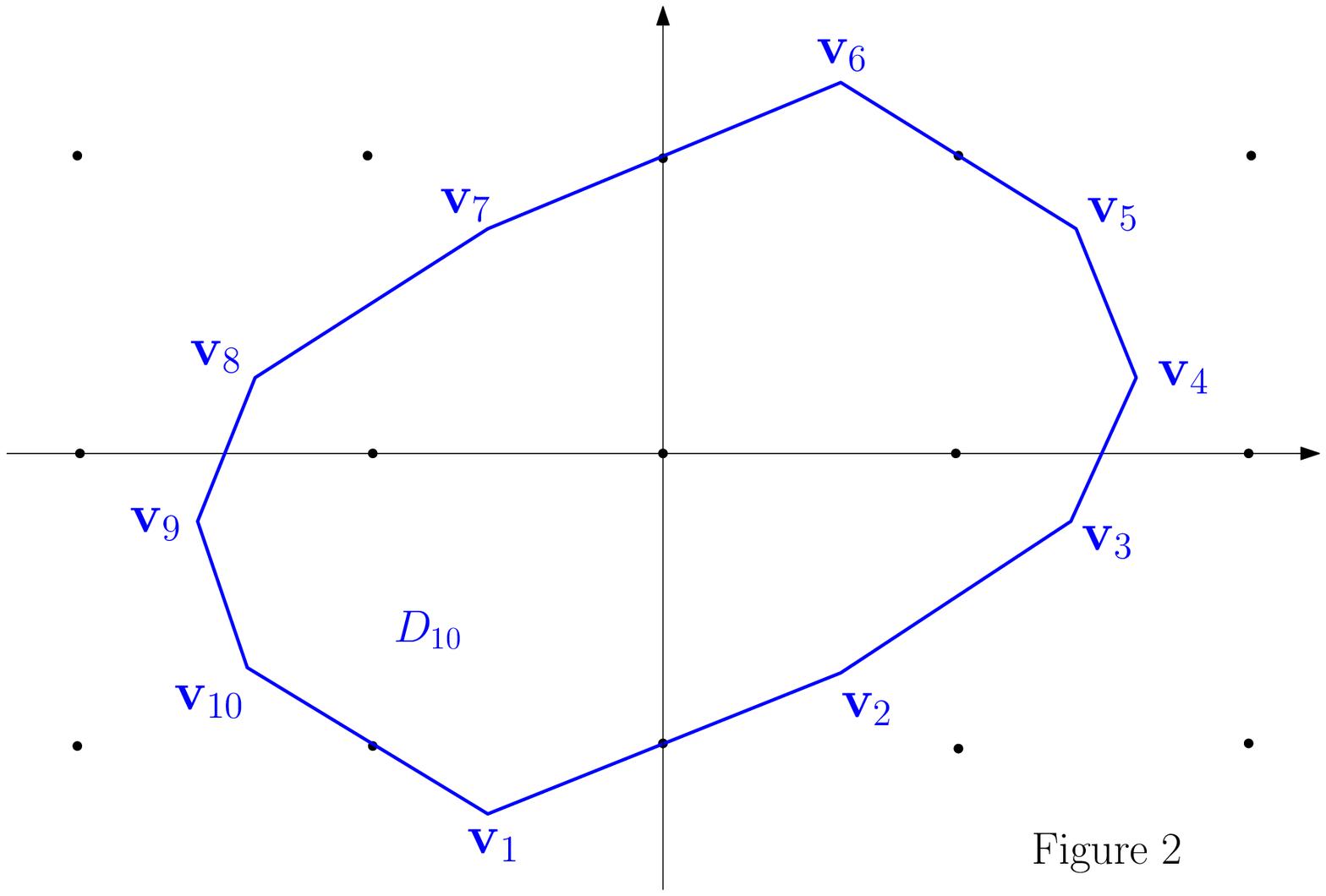}
\end{figure}

It can be easily verified that
$${\bf u}_i=\mbox{${1\over 2}$}({\bf v}_i+{\bf v}_{i+1})\in \mbox{${1\over 2}$}\Lambda, \quad i=1, 2, \ldots, 10, $$
where ${\bf v}_{11}={\bf v}_1$, and
$${\rm vol}(D_{10})=5.$$
It follows from Lemma 2 that $D_{10}+\Lambda$ is a five-fold lattice tiling. Combined with (2), it can be deduced that
$$\tau^*(D_{10})=5.$$
Theorem 1 is proved.\hfill{$\Box$}

\vspace{0.6cm}
\noindent
{\Large\bf 3. Some Comparisons and Generalizations}

\bigskip
\noindent
It is interesting to make comparisons with multiple packings and multiple coverings. Let $O$ denote the unit circular disk. It was discovered by Blunden \cite{blun57, blun63} that
$$\delta_k(O)=k\cdot \delta (O)$$
is no longer true when $k\ge 5$, and
$$\theta_k(O)=k\cdot \theta (O)$$
is no longer true when $k\ge 3$. So, the packing case is rather similar to the tilings, while the covering case is much different!

On the other hand, for every two-dimensional convex domain $D$ it was proved by Cohn \cite{cohn76}, Bolle \cite{bolle89} and Groemer \cite{groem86} that
$$\lim_{k\to\infty }{{\delta_k(D)}\over k}=\lim_{k\to\infty }{{\theta_k(D)}\over
k}=1.$$
In other words, from the density point of view, when the multiplicity is big there is no much difference among packing, covering and tiling!

\bigskip
Let $P$ denote an $n$-dimensional centrally symmetric convex polytope with centrally symmetric facets and let $V$ denote a $(n-2)$-dimensional face of $P$. We call the collection of all those facets of $P$ which contain a translate of $V$ as a subface a belt of $P$.

\medskip
In 1980, P. McMullen \cite{mcmu} proved the following criterion for parallelohedra.

\medskip\noindent
{\bf Lemma 3.} {\it  A convex body $K$ is a parallelohedron if and only if it is a centrally symmetric polytope with centrally symmetric facets and each belt contains four or six facets.}

\medskip\noindent
{\bf Proof of Theorem 2.} For convenience, we write $\mathbb{E}^n=\mathbb{E}^2\times \mathbb{E}^{n-2}$. Let $P_{2m}$ be a centrally symmetric $2m$-gon ($m\ge 4$) such that $P_{2m}+\mathbb{Z}^2$ is a $k$-fold lattice tiling of $\mathbb{E}^2$, let $I^{n-2}$ denote the unit cube $\{ (x_3, x_4, \ldots, x_n):\ |x_i|\le {1\over 2}\}$ in $\mathbb{E}^{n-2}$, and define
$$P=P_{2m}\times I^{n-2}.$$
It is easy to see that $P+\mathbb{Z}^n$ is a $k$-fold lattice tiling of $\mathbb{E}^n$.

Let ${\bf v}_1$, ${\bf v}_2$, $\ldots$, ${\bf v}_{2m}$ be the $2m$ vertices of $P_{2m}$ and let $G_1$, $G_2$, $\ldots $, $G_{2m}$ denote the $2m$ edges of $P_{2m}$, and define
$$V={\bf v}_1\times I^{n-2}$$
and
$$F_i=G_i\times I^{n-2}.$$
Clearly, $\{ F_1, F_2, \ldots , F_{2m}\}$ is a belt of $P$ with $2m$ facets. Therefore, by McMullen's criterion it follows that $P$ is not a parallelohedron in $\mathbb{E}^n$. In particular, the octagon $D_8$ and the decagon $D_{10}$ defined in the proof of Theorem 1 produce non-parallelohedral five-fold lattice tiles $D_8\times I^{n-2}$ and $D_{10}\times I^{n-2}$, respectively.

Theorem 2 is proved. \hfill{$\Box$}

\vspace{0.6cm}\noindent
{\bf Acknowledgements.} The authors are grateful to Professor Weiping Zhang for calling their attention to Rao's announcement, to Professors Mihalis Kolountzakis, Sinai Robins and G\"unter M. Ziegler for their comments, and to the referee for his revision suggestions. This work is supported by 973 Program 2013CB834201.

\bibliographystyle{amsplain}

\vspace{0.8cm}
\noindent
Qi Yang, School of Mathematical Sciences, Peking University, Beijing 100871, China

\vspace{0.3cm}
\noindent
Corresponding author:

\noindent
Chuanming Zong, Center for Applied Mathematics, Tianjin University, Tianjin 300072, China

\noindent
Email: cmzong@math.pku.edu.cn

\end{document}